\documentclass[12pt]{article}
\usepackage{amsmath,amssymb,amsfonts,amsthm}
\newenvironment{myabstract}{\par\noindent
{\bf Abstract . } \small }
{\par\vskip8pt minus3pt\rm}
\newcounter{item}[section]
\newcounter{kirshr}
\newcounter{kirsha}
\newcounter{kirshb}
\newenvironment{enumroman}{\setcounter{kirshr}{1}
\begin{list}{(\roman{kirshr})}{\usecounter{kirshr}} }{\end{list}}
\newtheorem{theorem}{Theorem}[section]

\newtheorem{lemma}[theorem]{Lemma}

\newenvironment{demo}[1]{\noindent{\bf #1.}\upshape\mdseries}
{\nopagebreak{\hfill\rule{2mm}{2mm}\nopagebreak}\par\normalfont}
\theoremstyle{definition}

\newtheorem{definition}[theorem]{Definition}

\def\Fr{{\mathfrak{Fr}}}
\def\Sg{{\mathfrak{Sg}}}

\def\K{{\mathfrak{K}}}

\def\K{{\bf K}}

\def\(R)RA{{\bf (R)RA}}

\def\A{{\mathfrak{A}}}
\def\B{{\mathfrak{B}}}

\def\D{{\mathfrak{D}}}

\def\Fr{{\mathfrak{Fr}}}
\def\Sg{{\mathfrak{Sg}}}

\def\K{{\mathfrak{K}}}

\def\K{{\bf K}}

\def\(R)RA{{\bf (R)RA}}

\def\A{{\mathfrak{A}}}
\def\B{{\mathfrak{B}}}

\def\D{{\mathfrak{D}}}
\def\Ig{{\mathfrak{Ig}}}

\def\E{{\mathfrak{E}}}

\title{Varieties of algebras without the amalgamation property} 
\author{Tarek Sayed Ahmed\\
Department of Mathematics, Faculty of Science,\\ 
Cairo University, Giza, Egypt.
  }
\begin{document}
\maketitle
\begin{myabstract}   Let $\alpha$ be an ordinal and $\kappa$ be a cardinal, both infinite, such that $\kappa\leq |\alpha|.$
For $\tau\in {}^{\alpha}\alpha$, let $sup(\tau)=\{i\in \alpha: \tau(i)\neq i\}$. Let $G_{\kappa}=\{\tau\in {}^{\alpha}\alpha: |sup(\tau)|< \kappa\}$.
We consider variants of polyadic equality algebras by taking cylindrifications on $\Gamma\subseteq \alpha$, $|\Gamma|<\kappa$ and substitutions restricted 
to $G_{\kappa}$. Such algebras are also enriched with generalized diagonal elements.
We show that for any variety $V$ containing the class of representable algebas and 
satisfying a finite schema of equations, $V$ fails to have the amalgamation
property. In particular, many varieties of Halmos' quasi-polyadic equality algebras and Lucas' extended cylindric algebras 
(including that of the representable algebras) fail to have the amalgamation property. 

\end{myabstract}

The most generic examples of algebraisations of first order logic are Tarski's cylindric algebras and Halmos' polyadic algebras.
Both algebras are well known and widely used.
Polyadic algebras were introduced by Halmos \cite{Halmos} to provide an algebraic reflection
of the study of first order logic without equality. Later the algebras were enriched by 
diagonal elements to permit the discussion of equality. 
That the notion is indeed an adequate reflection of first order logic was 
demonstrated by Halmos' representation theorem for locally finite polyadic algebras 
(with and without equality). Tarski proved an analogous result for locally finite cylindric algebras.
Daigneault and Monk  
proved a strong extension of Halmos' theorem, namely,  
every polyadic algebra of infinite dimension (without equality) is representable \cite{DM}. However, not every cylindric algebra
is representable. In fact, the class of infinite dimensional representable algebras is not axiomatizable by any finite schema, a classical result of Monk.
This is a point (among others) where the two theories deviate. Monk's result was considerably strengthened by Andr\'eka by showing that
there is an inevitable degree of complexity in any 
axiomatization of the class of representable cylindric algebras.
In particular, any universal axiomatization of the class of representable quasipolyadic algebras must contain infinitely many variables.
The representation theorem of Diagneualt and Monk - a typical Stone-like representation theorem - shows 
that the notion of polyadic algebra is indeed an adequate 
reflection of Keisler's predicate logic ($KL$).
$KL$ is a proper extension of first order logic without equality,
obtained when the bound on the number of variables in formulas is relaxed; 
and  accordingly allowing the following as extra operations on formulas:
Quantification on infinitely many variables 
and simultaneous substitution of (infinitely many) 
variables for variables. 
Adding equality to $KL$, proved problematic as illustrated 
algebraically by Johnson \cite{J69}.  
In op.cit, Johnson showed that the class of 
representable polyadic algebras with equality is not closed under ultraproducts, hence 
this class is not elementary,  i.e. cannot be axiomatized by any set of first order
sentences.  However one can still hope for a nice axiomatization of the {\it variety generated}
by the class of polyadic equality algebras. A subtle recent (negative) 
result in this direction is 
N\'emeti - S\'agi's \cite{NS}:
In sharp contrast to $KL$, the validities of $KL$ {\it with equality} 
cannot be recaptured by any set of schemas analogous to Halmos' schemas, 
let alone a finite one. 
In particular, the variety generated by the class of representable polyadic 
algebras with equality cannot be axiomatized by a finite schema of equations. 
The latter answers a question originally raised by Craig \cite{C}.

It is interesting (and indeed natural) to ask for algebraic versions 
of model theoretic results, 
other than completeness. 
Examples include interpolation theorems and omitting types theorems. 
Unlike the cylindric case, omitting types for polyadic algebras prove problematic. This is the case 
because polyadic algebras
of infinite dimension 
have uncountably many operations, and omitting types arguments-  Baire Category arguments at heart - are very much tied to countability.
On the other hand, Daigneault succeeded in stating and proving versions of Beth's and Craig's 
theorems. This was done by proving the algebraic analogue of Robinson's joint 
consistency  theorem: {\it Locally finite} polyadic algebras (with and without equality) 
have  the amalgamation property. 
Later Johnson removed the condition of local finiteness, proving that polyadic algebras 
{\it without} equality have the {\it strong} amalgamation property \cite{J70} . 
With this stronger result, Robinson's, Beth's and Craig's theorems hold for $KL$. 

Yet another point where the two theories deviate, Pigozzi \cite{P} proves that the class of representable cylindric algebras fails to have the amalgamation property.
This shows that certain infinitary 
algebraisable extensions of first order logic, 
the so-called
typless logics (or finitary logics with infinitary relations) fail to have the interpolation property.  
Further negative results concerning various amalgamation
properties for cylindric-like algebras of relations can be found in \cite{AUU}, \cite{MSone}, \cite{MStwo}.

Motivated by the quest for algebraisations that posses the positive properties of both polyadic algebras and cylindric algebras, in this paper we show, 
using basically Pigozzi's techniques appropriately modified,  
that the interpolation property fails for many variants of $KL$ {\it with} equality, contrasting the equality free case \cite{AU}.
In such variants, formulas of infinite length are allowed, but quantification and substitutions are only allowed for $<\kappa$ many variables where
$\kappa$ is a fixed beforehand infinite cardinal.  
Also (generalized) equality is available. Such logics are (natural) extensions of the typeless logics corresponding to cylindric algebras.

Our proof is algebraic adressing the amalgamation property for certain variants of the class of polyadic equality algebras, that are also proper expansions of cylindric algebras.
From our proof it can be easily destilled that many varieties of algebraic logics existing in the literature fail to have the amalgamation property.
Examples include Halmos' quasi-polyadic equality algebras and Lucas' extended cylindric algebras.  These results are new.

\section{Results and proofs}

Let $\alpha$ be an ordinal and $\kappa$ be a cardinal, both infinite,  
such that $\kappa\leq |\alpha|$. For $\tau\in {}^{\alpha}\alpha$, let $sup(\tau)=\{i\in \alpha: \tau(i)\neq i\}$.
Let $G_{\kappa}=\{\tau\in {}^{\alpha}\alpha: |sup(\tau)|<\kappa\}$. Clearly $G_k$ is a semigroup under the operation of composition; in fact it is a monoid. 
We write $\Gamma\subseteq_{\kappa}\alpha$
if $\Gamma\subseteq \alpha$ and $|\Gamma|<\kappa$. Let 
$N=\{\E\subseteq \alpha\times \alpha: \text { $\E$ is an equivalence relation on $\alpha$ and } |\{i<\alpha:i/\E\neq \{i\}\}|<\kappa\}.$

\begin{definition} By a $\kappa$ generalized polyadic equality algebra
dimension $\alpha$, or a $PEA_{\kappa, \alpha}$ for short,
we understand an algebra of the form
$$\A=\langle A,+,\cdot ,-,0,1,{\sf c}_{(\Gamma)},{\sf s}_{\tau}, {\sf d}_{\E} \rangle_{\Gamma\subseteq_{\kappa} \alpha ,\tau\in G_{\kappa}, \E\in N}$$
where ${\sf c}_{(\Gamma)}$ ($\Gamma\subseteq_{\kappa}  \alpha$) and ${\sf s}_{\tau}$ ($\tau\in G_{\kappa})$ are unary 
operations on $A$, ${\sf d}_{E}\in A$ ($E\in N$), such that postulates  
below hold for $x,y\in A$, $\tau,\sigma\in G,$
$\Gamma, \Delta\subseteq_{\kappa} \alpha$, $E\in N$, and all $i,j\in \alpha$.

\begin{enumerate}
\item $\langle A,+,\cdot ,-,0,1\rangle$ is a boolean algebra
\item ${\sf c}_{(\Gamma)}0=0$

\item $x\leq {\sf c}_{{(\Gamma)}}x$

\item ${\sf c}_{(\Gamma)}(x\cdot {\sf c}_{(\Gamma)}y)={\sf c}_{(\Gamma)}x\cdot {\sf c}_{(\Gamma)}y$

\item ${\sf c}_{(\Gamma)}{\sf c}_{(\Delta)}x={\sf c}_{(\Gamma\cup \Delta)}x$

\item ${\sf s}_{\tau}$ is a boolean endomorphism

\item ${\sf s}_{Id}x=x$

\item ${\sf s}_{\sigma\circ \tau}={\sf s}_{\sigma}\circ {\sf s}_{\tau}$

\item if $\sigma\upharpoonright (\alpha\sim \Gamma)=\tau\upharpoonright (\alpha\sim \Gamma)$, then 
${\sf s}_{\sigma}{\sf c}_{(\Gamma)}x={\sf s}_{\tau}{\sf c}_{(\Gamma)}x$


\item If $\tau^{-1}\Gamma=\Delta$ and $\tau\upharpoonright \Delta $ is
one to one, then ${\sf c}_{(\Gamma)}{\sf s}_{\tau}x={\sf s}_{\tau}{\sf c}_{(\Delta )}x$



\item ${\sf d}_{I}=1$ where $I=Id\upharpoonright \alpha\times \alpha$

\item ${\sf c}_{(\Gamma)}{\sf d}_E={\sf d}_F$ where $F=E\cap {}^2(\alpha\sim \Gamma)\cup Id\upharpoonright \alpha\times \alpha$

\item ${\sf s}_{\tau}{\sf d}_{E}={\sf d}_{F}$  where $F=\{(\tau(i), \tau(j)): (i,j)\in E\}\cup Id\upharpoonright \alpha\times \alpha$

\item $x\cdot {\sf d}_{ij}\leq {\sf s}_{[i|j]}x$ 


\end{enumerate}

\end{definition}
In the above definition, and elsewhere throughout the paper,  ${\sf d}_{ij}$ denotes the element ${\sf d}_E$ 
where $E$ is the equivalence relation relating $i$ to $j$, and everything else only to itself.
For a class $\K$ of algebras, $S\K$ stands for the class of all subalgebras of algebras in $\K$, $P\K$ is the class of products of algebras in $\K$
and $H\K$ is the class of all homomorphic images of algebras in $\K$. The class of representable algebras is defined via set - theoretic operations
on sets of $\alpha$-ary sequences. Let $U$ be a set. 
For $\Gamma\subseteq \alpha$,  $\tau\in {}^{\alpha}\alpha, i, j\in \alpha$  and $E\in N$,  we set
$${\sf c}_{(\Gamma)}X=\{s\in {}^{\alpha}U: \exists t\in X\ t(j)=s(j)  \ \ \forall \\ j\notin \Gamma\}.$$
$${\sf s}_{\tau}X=\{s\in {}^{\alpha}U: s\circ \tau\in X\}.$$
$${\sf d}_{ij}=\{s\in {}^{\alpha}U: s_i=s_j\}.$$
$${\sf d}_{E}=\{s\in {}^{\alpha}U:  s_i=s_j \ \ \forall (i,j)\in E\}.$$
Note that ${\sf d}_E=\bigcap_{(i,j)\in E}{\sf d}_{ij}$. For a set $X$, let $\B(X)$ be the boolean set algebra $(\wp(X), \cap, \cup, \sim).$
The class of representable $G_k$ polyadic equality algebras, or $RPEA_{\kappa, \alpha}$ is defined by
$$SP\{\langle \B(^{\alpha}U), {\sf c}_{(\Gamma)}, {\sf s}_{\tau}, {\sf d}_{E}\rangle: E\in N, \Gamma\subseteq_{\kappa} \alpha, \tau\in G_{\kappa},
\ \,  U \text { a set }\}.$$ 
We make the following observations:

\begin{itemize}

\item $RPEA_{\kappa, \alpha}\subseteq PEA_{\kappa, \alpha}$, and the inclusion is proper \cite{p}.

\item If $\A\in PEA_{\kappa, \alpha}$ then $\A$ has a cylindric reduct and indeed this reduct is a cylindric algebra of dimension $\alpha$. 
In fact, $\A$ has a quasipolyadic equality reduct obtained by restricting the operations to 
finite quantifiers (cylindrifications) , finite substitutions and ordinary diagonal elements, i.e. the ${\sf d}_{ij}$'s.

\item if $G_{\kappa}$ contains one infinitary substitution then $RPEA_{\kappa, \alpha}$ is not closed under ultraproducts \cite{Sain}, 
hence is not closed under $H$, lest it be a variety.

\end{itemize}
For what follows, we need:

\begin{definition}
Let $K\subseteq V$ be classes of algebras.
$K$ is said to have the 
amalgamation property, or
$AP$ for short, with respect to $V$, if for all ${\A}_0$, ${\A}_1$ 
and ${\A}_2\in K$, and all
monomorphisms $i_1$ and $i_2$ of ${\A}_0$ into ${\A}_1$,
${\A}_2$,  respectively, there exists ${\A}\in V$, a
monomorphism $m_1$ from ${\A}_1$ into $\A$ and a
monomorphism $m_2$ from ${\A}_2$ into $\A$ such that
$m_1\circ i_1=m_2\circ i_2$. 
\end{definition}
We will show that for any variety $\K$, $RPEA_{\kappa, \alpha}\subseteq \K\subseteq PEA_{\kappa, \alpha}$, 
$\K$ fails to have the amalgamation property with respect to $PEA_{\kappa, \alpha}$.  
For motivations of studying such algebras, and similar reducts of polyadic equality algebras, initiated by Craig \cite{C},  see \cite{Tarskian}, \cite{AU}, \cite{Bulletin}, \cite{Sain}, \cite{SG}.
Amalgamation in varieties can be pinned down to congruences on free algebras.
Congruences correspond to ideals. This prompts:

\begin{definition} Let $\A\in PEA_{\kappa, \alpha}$. A subset $I$ of $\A$
in an ideal if the following conditions are satisfied:

\begin{enumroman}

\item $0\in I,$

\item If $x,y\in I$, then $x+y\in I,$

\item If $x\in I$ and $y\leq x$ then $y\in I,$

\item For all $\Gamma\subseteq_{\kappa} \alpha$ and $\tau\in G_{\kappa}$ if $x\in I$
then ${\sf c}_{(\Gamma)}x$ and ${\sf s}_{\tau}x\in I$.

\end{enumroman}
\end{definition}
It can be checked that ideals function properly, that is ideals correspond to congruences 
the usual way.
For $X\subseteq \A$, the ideal generated by $X$, $\Ig^{\A}X$ is the smallest ideal containing 
$X$, i.e the intersection of all ideals containing $X$. We let $\Sg^{\A}X$ and sometimes ${\A}^{(X)}$ denote the subalgebra of 
$\A$ generated by $X$. 

\begin{lemma} Let $\A\in PEA_{\kappa, \alpha}$ and $X\subseteq A$.
Then 
$\Ig^{\A}X=\{y\in A: y\leq {\sf c}_{(\Gamma)}(x_0+\ldots x_{k-1})\}: 
\text{ for some  } k\in \omega,  x\in {}^kX \text { and }\Gamma\subseteq_{\kappa} \alpha\}.$
\end{lemma}
\begin{demo}{Proof}
Let $H$ denote the set of elements on the right hand side.
It is easy to check $H\subseteq \Ig^{\A}X$.
Conversely, assume that $y\in H,$ $\Gamma\subseteq_{\kappa} \alpha.$ 
It is clear that ${\sf c}_{(\Gamma)}y\in H$.
$H$ is closed under substitutions, since for any $\tau\in G_{\kappa}$, any $x\in A$ there exists $\Gamma\subseteq_{\kappa} \alpha$
such that ${\sf s}_{\tau}x\leq {\sf c}_{(\Gamma)}x$. Indeed $sup(\tau)$ is such a $\Gamma$.
Now let $z, y\in H$. Assume that $z\leq {\sf c}_{(\Gamma)}(x_0+\ldots x_{k-1})$
and $y\leq {\sf c}_{(\Delta)}(y_0+\ldots y_{l-1}),$
then $$z+y\leq {\sf c}_{(\Gamma\cup \Delta)}(x_0+\ldots x_{k-1}+ y_0\ldots +y_{l-1}).$$
The Lemma is proved. 
\end{demo}
Fixing $\alpha$ and $\kappa$ throughout, in what follows we denote $(R)PEA_{\kappa, \alpha}$ simply by $(R)PEA.$
The following about ideals  will be frequently used.
\begin{itemize}
\item If $\A\subseteq \B$ are $PEA$'s and $I$ is an ideal of $\A$, then $\Ig^{\B}(I)=\{b\in B: \exists a\in I  (b\leq a)\}.$
\item If $I$ and $J$ are ideals of a $PEA$ then the ideal generated by $I\cup J$ is
$\{x: x\leq i+j \text { for } i\in I, j\in J\}.$
\end{itemize}

For a class $\K$ and a set $X$, $\Fr_X\K$ denotes the $\K$ algebra 
freely generated by $X$, or the $\K$ free algebra on $|X|$ generators.
As a wide spread custom, we identify $X$ with $|X|$.
We understand the notion of free algebras in the sence of \cite{HMT1} Definition 0.4.19. In particular, free $\K$ algebras
may not be in $\K$. However, they are always in $HSP(\K)$, the variety generated by $\K$. 
We write $R\in Co\A$ if $R$ is a congruence relation $\A$. 
For $X\subseteq A$, then by $(\A/R)^{(X)}$ we undertand the subalgebra of 
$\A/R$ generated by $\{x/R: x\in X\}.$
Since our algebras have cylindric reducts, in what follows we use freely results of Henkin Monk and Tarski's treatise \cite{HMT1} on the arithmetic of cylindric algebras.
We now formulate and prove our main result:
\begin{theorem} Let $\K$ be a variety such that $RPEA\subseteq \K\subseteq PEA$. Then $\K$ 
does not have $AP$ with respect to $PEA.$
\end{theorem}

\begin{demo}{Proof} The proof is an adaptation of Pigozzi's techniques  for showing failure of the amalgamation property for cylindric algebras
\cite{P}. Seeking a contradiction assume that $\K$ has $AP$ with respect
to $PEA$. Let $\A=\Fr_4PEA.$ 
Let $r, s$ and $t$ be defined as follows:
$$ r = {\sf c}_0(x\cdot {\sf c}_1y)\cdot {\sf c}_0(x\cdot -{\sf c}_1y),$$
$$ s = {\sf c}_0{\sf c}_1({\sf c}_1z\cdot {\sf s}^0_1{\sf c}_1z\cdot -{\sf d}_{01}) + {\sf c}_0(x\cdot -{\sf c}_1z),$$
$$ t = {\sf c}_0{\sf c}_1({\sf c}_1w\cdot {\sf s}^0_1{\sf c}_1w\cdot -{\sf d}_{01}) + {\sf c}_0(x\cdot -{\sf c}_1w),$$
where $ x, y, z, \text { and } w$ are the first four free generators
of $\A$. 
Then $r\leq s\cdot t$. This  inequality is proved by Pigozzi, whose proof we include. Indeed put
$$a = x\cdot {\sf c}_1 y\cdot  -{\sf c}_0(x\cdot - {\sf c}_1z),$$
$$b = x\cdot  - {\sf c}_1 y\cdot -{\sf c}_0(x\cdot - {\sf c}_1z).$$
Then we have
\begin{equation*}
\begin{split}
 {\sf c}_1 a\cdot  {\sf c}_1 b & \leq {\sf c}_1 (x \cdot {\sf c}_1y)\cdot  {\sf c}_1(x\cdot -{\sf c}_1y)
  \,\,\ \textrm{by \cite{HMT1}1.2.7}\\
&  = {\sf c}_1 x\cdot {\sf c}_1 y\cdot  {\sf c}_1 x\cdot -{\sf c}_1y \,\,\,\, \textrm{by \cite{HMT1} 1.2.11}
\end{split}
\end{equation*}
and so
\begin{equation}\label{p1}
\begin{split}
{\sf c}_1 a\cdot  {\sf c}_1 b = 0.
\end{split}
\end{equation}
From the inclusion $x\cdot -{\sf c}_1 z \leq {\sf c}_0 (x\cdot -{\sf c}_1z)$ we get
$$x\cdot -{\sf c}_0 (x\cdot -{\sf c}_1z) \leq {\sf c}_1 z.$$
Thus $a, b \leq {\sf c}_1z$ and hence, by \cite{HMT1} 1.2.9,
\begin{equation}\label{p2}
\begin{split}
{\sf c}_1 a, {\sf c}_1 b \leq {\sf c}_1z.
\end{split}
\end{equation}
 We now compute:
\begin{equation*}
\begin{split}
{\sf c}_0 a\cdot  {\sf c}_0 b & \leq {\sf c}_0 {\sf c}_1 a \cdot 
{\sf c}_0 {\sf c}_1 b \,\,\,\,\,\  \textrm{by \cite{HMT1}
 1.2.7} \\
& = {\sf c}_0 {\sf c}_1 a \cdot  {\sf c}_1 {\sf s}^0_1 {\sf c}_1 b \,\,\,\,\,\  \textrm{by \cite{HMT1} 1.5.8
(i),\,\,\,\,\,\ \cite{HMT1} 1.5.9 (i)}\\
& = {\sf c}_1({\sf c}_0{\sf c}_1 a\cdot  {\sf s}^0_1 {\sf c}_1 b)\\
& = {\sf c}_0{\sf c}_1({\sf c}_1a\cdot {\sf s}^0_1{\sf c}_1b)\\
& = {\sf c}_0{\sf c}_1[{\sf c}_1a\cdot {\sf s}^0_1{\sf c}_1b\cdot (-{\sf d}_{01} + {\sf d}_{01})\\
& = {\sf c}_0{\sf c}_1[({\sf c}_1a\cdot {\sf s}^0_1{\sf c}_1b\cdot -{\sf d}_{01}) + ({\sf c}_1a\cdot {\sf s}^0_1{\sf c}_1b.{\sf d}_{01})]\\
& = {\sf c}_0{\sf c}_1[({\sf c}_1a\cdot {\sf s}^0_1{\sf c}_1b\cdot -{\sf d}_{01}) + ({\sf c}_1a\cdot {\sf c}_1b\cdot {\sf d}_{01})] \,\,\,\,\
\textrm{by \cite{HMT1}
1.5.5}\\
& = {\sf c}_0{\sf c}_1 ({\sf c}_1a\cdot {\sf s}^0_1{\sf c}_1b\cdot -{\sf d}_{01}) \,\,\,\,\  \textrm{by (\ref{p1})}\\
& \leq {\sf c}_0{\sf c}_1 ({\sf c}_1z\cdot {\sf s}^0_1{\sf c}_1z\cdot -{\sf d}_{01}) \,\,\,\,\  \textrm{by (\ref{p2}), \cite{HMT1} 1.2.7}\\
\end{split}
\end{equation*}
We have proved that
$$ {\sf c}_0[ x\cdot  {\sf c}_1y\cdot -{\sf c}_0(x\cdot -{\sf c}_1z)]\cdot {\sf c}_0[x\cdot -{\sf c}_1y\cdot -{\sf c}_0(x\cdot -{\sf c}_1z)] \leq
{\sf c}_0{\sf c}_1({\sf c}_1z\cdot {\sf s}^0_1{\sf c}_1z\cdot -{\sf d}_{01}).$$
In view of \cite{HMT1} 1.2.11 this gives
$$ {\sf c}_0( x\cdot  {\sf c}_1y)\cdot {\sf c}_0(x\cdot -{\sf c}_1y)\cdot -{\sf c}_0(x\cdot -{\sf c}_1z) \leq
{\sf c}_0{\sf c}_1({\sf c}_1z\cdot {\sf s}^0_1{\sf c}_1z\cdot -{\sf d}_{01}).$$
The conclusion now follows.
Let $X_1=\{x, y\}$
and $X_2=\{x,z, w\}.$ Then  
\begin{equation}\label{p15}
\begin{split}
\A^{(X_1\cap X_2)}=\Sg^{\A}\{x\}.
\end{split}
\end{equation}
We have 
\begin{equation}
\begin{split}
r\in A^{(X_1)} \text { and }s,t\in A^{(X_2)}.
\end{split}
\end{equation}
Let $R$ be an ideal of $\A$ such that
\begin{equation}\label{p17}
\begin{split}
\A/R\cong \Fr_4\K_{\alpha}.
\end{split}
\end{equation}
Since $r\leq s\cdot t$ we have  
\begin{equation}\label{p18}
\begin{split}
r\in \Ig^{\A}\{s\cdot t\}\cap A^{(X_1)}.
\end{split}
\end{equation}
Let
\begin{equation}\label{p19}
\begin{split}
M=\Ig^{\A^{(X_2)}}[\{s\cdot t\}\cup (R\cap A^{(X_2)})];
\end{split}
\end{equation}
\begin{equation}\label{p20}
\begin{split}
N=\Ig^{\A^{(X_1)}}[(M\cap A^{(X_1\cap X_2)})\cup (R\cap A^{(X_1)})].
\end{split}
\end{equation}
Then we have
\begin{equation}\label{p21}
\begin{split}
R\cap A^{(X_2)}\subseteq M\text { and }R\cap A^{(X_1)}\subseteq N.
\end{split}
\end{equation}
From the first of these inclusions we get
$$M\cap A^{(X_1\cap X_2)}\supseteq (R\cap A^{(X_2)})\cap A^{(X_1\cap X_2)}=(R\cap A^{(X_1)})
\cap A^{(X_1\cap X_2)}.$$
By $(8)$ we have 
$$N\cap A^{(X_1\cap X_2)}=M\cap A^{(X_1\cap X_2)}.$$

For $R$ an ideal of $\A$ and $X\subseteq A$, by $(\A/R)^{(X)}$ we understand the subalgebra of 
$\A/R$ generated by $\{x/R: x\in X\}.$  
Define $$\theta: \A^{(X_1\cap X_2)}\to \A^{(X_1)}/N$$
by $$a\mapsto a/N.$$
Then $ker\theta=N\cap A^{(X_1\cap X_2)}$ and $Im\theta=(\A^{(X_1)}/N)^{(X_1\cap X_2)}$.
It follows that $$\bar{\theta}:\A^{(X_1\cap X_2)}/N\cap {}A^{(X_1\cap X_2)}\to (\A^{(X_1)}/N)^{(X_1\cap X_2)}$$
defined by
$$a/N\cap {}A^{X_1\cap X_2)}\mapsto a/N$$
is a well defined isomorphism.
Similarly
$$\bar{\psi}:\A^{(X_1\cap X_2)}/M\cap {}A^{(X_1\cap X_2)}\to (\A^{(X_2)}/M)^{(X_1\cap X_2)}$$
defined by
$$a/M\cap {}A^{X_1\cap X_2)}\mapsto a/M$$
is also a well defined isomorphism.
But $$N\cap A^{(X_1\cap X_2)}=M\cap A^{(X_1\cap X_2)},$$
Hence
$$\phi: (\A^{(X_1)}/N)^{(X_1\cap X_2)}\to (\A^{(X_2)}/M)^{(X_1\cap X_2)}$$
defined by
$$a/N\mapsto a/M$$
is a well defined isomorphism. 
Now
$(\A^{(X_1)}/N)^{(X_1\cap X_2)}$ embeds into ${\A}^{(X_1)}/N$ via the inclusion map; it also embeds in $\A^{(X_2)}/M$ via $i\circ \phi$ where $i$ 
is also the inclusion map.
For brevity let $\A_0=(\A^{(X_1)}/N)^{(X_1\cap X_2)}$, $\A_1={\A}^{(X_1)}/N$ and $\A_2={\A}^{(X_2)}/M$ and $j=i\circ \phi$.
Then $\A_0$ embeds in $\A_1$ and $\A_2$ via $i$ and $j$ respectively. Now observe that $\A_1$, $\A_2$ and $\A_0$ are in $\K$.
So by assumption, there exists an amalgam, i.e
there exists  $\B\in PEA$ and monomorphisms $f$ and $g$ from $\A_1$ and $\A_2$ respectively to 
$\B$ such that
$f\circ i=g\circ j$.
Let $$\bar{f}:\A^{(X_1)}\to \B$$ be defined by $$a\mapsto f(a/N)$$ and $$\bar{g}:\A^{(X_2)}\to \B$$ 
be defined by $$a\mapsto g(a/M).$$
Let $\B'$ be the algebra generated by $Imf\cup Im g$.
Then $\bar{f}\cup \bar{g}\upharpoonright X_1\cup X_2\to \B'$ is a function since $\bar{f}$ and $\bar{g}$ coincide on $X_1\cap X_2$.
By freeness of $\A$, there exists $h:\A\to \B'$ such that $h\upharpoonright_{X_1\cup X_2}=\bar{f}\cup \bar{g}$.
Let $P=kerh $. Then it is not hard to check that 
\begin{equation}\label{p23}
\begin{split}
P\cap A^{(X_1)}=N,
\end{split}
\end{equation}
and 
\begin{equation}\label{p24}
\begin{split}
P\cap A^{(X_2)}=M.
\end{split}
\end{equation}

In view of $(4),(7),(11)$ we have $s\cdot t\in P$ and hence by $(6)$ $r\in P$. 
Consequently from $(4)$ and $(11)$  we get $r\in N$. From $(8)$ 
there exist elements
\begin{equation}\label{p25}
\begin{split}
u\in M\cap A^{(X_1\cap X_2)}
\end{split}
\end{equation}
and $b\in R$ such that
\begin{equation}\label{p26}
\begin{split}
r\leq u+b.
\end{split}
\end{equation}
Since $u\in M$ by $(7)$ there is a $\Gamma\subseteq_{\kappa} \alpha$ and 
$c\in R$ such that
$$u\leq {\sf c}_{(\Gamma)}(s\cdot t)+c.$$
Let $\{x', y', z', w'\}$ be the first four generators of $\D=\Fr_{4}\K$. Let $h$ be the homomorphism from $\A$ to $\D$ be such that
$h(i)=i'$ for $i\in \{x,y,w,z\}$.
Notice that $ker h=R$. Then $h(b)=h(c)=0$.
It follows that $$h(r)\leq h(u) \leq {\sf c}_{(\Gamma)}(h(s).h(t)).$$
Let $r'=h(r)$, $u'=h(u)$, $s'=h(s)$ and $t'=h(t).$ Let $$\B = ( \wp (^{\alpha}{\alpha}), \cup, \cap, \sim, \emptyset,
{^{\alpha}{\alpha}}, {\sf c}_{(\Gamma)}, {\sf s}_{\tau}, {\sf d}_\E)_{
 \Gamma\subseteq_{\kappa} \alpha, \tau\in  G_{\kappa}, \E\in N}$$ 
that is $\B$ is the  full set algebra in the
space $ {^{\alpha}{\alpha}}$. Let $E$ be the set of all equivalence
relations on $\alpha$, and for each $ R \in E $ set
$$ X_R = \{ \varphi : \varphi  \in {^{\alpha}{\alpha}} \,\,\,\
\textrm{and for all} \,\,\,\  \xi, \eta < \alpha , \varphi_\xi =
\varphi_\eta \,\,\,\ \textrm{iff} \,\,\,\  \xi R \eta \}.$$
More succintly
$$X_R=\{ \varphi \in {^{\alpha}{\alpha}}: ker \varphi=R\}.$$
Let
$$ C = \{ \bigcup_{R \in L} X_R : L \subseteq E \}. $$
$C$ is clearly closed under the formation of arbitrary unions, and
since
$$ \sim \bigcup_{R \in L} X_L =  \bigcup_{R \in E \sim L} X_R$$
for every $ L \subseteq E$, we see that $C$ is closed under the
formation of complements with respect to ${^{\alpha}{\alpha}}$. Thus
$ C $ is a Boolean subuniverse (indeed, a complete Boolean
subuniverse) of $\B$; moreover, it is obvious that
\begin{equation}\label{p4}
\begin{split}
X_R \,\,\,\ \textrm{is an atom of} \,\,\,\ ( C, \cup, \cap, \sim,
0, {^{\alpha}{\alpha}}) \,\,\,\,\ \textrm{for each} \,\,\,\ R \in E.
\end{split}
\end{equation}
For all $\E\in N$ we have ${\sf  d}_{\E} =
\bigcup \{ X_R : \E\subseteq R\in E \} $ and hence $
{\sf d}_{\E} \in C$. Also,
$$ {\sf c}_{(\Gamma)}X_R = \bigcup \{ X_S : S \in E, {}^{2}(\alpha \sim \Gamma) \cap S = {}^{2}(\alpha \sim \Gamma) \cap R \}$$
for any $ \Gamma\subseteq_{\kappa} \alpha$ and $ R \in E$. Thus, because ${\sf c}_{(\Gamma)}$
is completely additive, $C$ is closed under the operation $ {\sf c}_{(\Gamma)}$ 
for every $ \Gamma \subseteq_{\kappa} \alpha$. It is easy to show that $C$ is closed under substitutions.
For any $\tau\in G_{\kappa}$,
$${\sf s}_{\tau}X_R=\bigcup\{X_S: S\in E, \forall i, j<\alpha(iRj \Longleftrightarrow \tau(i)S\tau(j)\}.$$
The set on the right may of course be empty.
Since ${\sf s}_{\tau}$ is also completely additive, therefore, we have shown
that
\begin{equation}\label{p5}
\begin{split}
C \,\,\,\ \textrm{is a subuniverse of} \,\,\,\ \B .
\end{split}
\end{equation}
We now show that there is a subset $Y$ of $
{^{\alpha}{\alpha}}$ such that
\begin{equation}\label{p6}
\begin{split}
&X_{Id} \cap f(r') \neq 0 \,\,\ \textrm{for
every} \,\,\ f \in Hom (\D, \B) \\
& \textrm{such that} \,\,\ f(x') = X_{Id} \,\,\
\textrm{and} \,\,\ f(y') = Y,
\end{split}
\end{equation}
and also that for every $ \Gamma \subseteq_{\kappa} \alpha$, there are subsets $Z, W$ of $ {^{\alpha}{\alpha}}$ such that
\begin{equation}\label{p7}
\begin{split}
& X_{Id} \sim {\sf c}_{(\Gamma)}g(s'\cdot t') \neq 0 \,\,\
\textrm{for every} \,\,\ g\in Hom (\D, \B)\\
 & \textrm{such that}
\,\,\ g(x') = X_{Id}, g(z') = Z \,\,\ \textrm{and}
\,\,\ g(w') = W.
\end{split}
\end{equation}
Here $Hom(\D, \B)$ stands for the set of all homomorphisms from $\D$ to $\B$.
Let $\sigma\in {}^{\alpha}{\alpha}$ be such that 
$\sigma_0 = 0$, and $\sigma_\kappa = \kappa + 1$ for every non-zero
$\kappa < \omega$ and $\sigma j=j$ otherwise.
Let $ \tau = \sigma\upharpoonright (\alpha \sim \{0\}) \cup \{(0, 1)\}$. 
Then $
\sigma, \tau \in X_{Id}$. Take
$$ Y = \{\sigma\}.$$
Then
$$ \sigma \in X_{Id} \cap {\sf c}_1 Y \,\,\ \textrm{and}
\,\,\ \tau \in X_{Id} \sim {\sf c}_1 Y$$ and hence
\begin{equation}\label{p9}
\begin{split}
\sigma \in {\sf c}_0 (X_{Id} \cap {\sf c}_1 Y) \cap {\sf c}_0
(X_{Id} \sim {\sf c}_1 Y).
\end{split}
\end{equation}
Therefore, we have $ \sigma \in f(r')$ for every $ f \in Hom(\D, \B)$
such that $ f(x') = X_{Id}$ and $f(y') =Y$, and that
(\ref{p6}) holds.
We now want to show that for any given $ \Gamma \subseteq_{\kappa}
\alpha$ , there exist sets $ Z, W \subseteq {^{\alpha}{\alpha}}$ such
that (\ref{p7}) holds; it is clear that no generality is lost if we
assume that $0, 1 \in \Gamma$, so we make this assumption. Take
$$ Z = \{ \varphi : \varphi \in X_{Id},
\varphi_0 < \varphi_1 \} \cap {\sf c}_{(\Gamma)} \{Id\}$$ and
$$ W = \{ \varphi : \varphi \in X_{Id},
\varphi_0 > \varphi_1 \} \cap {\sf c}_{(\Gamma)} \{Id\}.$$
We show that
\begin{equation}\label{p10}
\begin{split}
Id \in X_{Id} \sim {\sf c}_{(\Gamma)} g(s'\cdot t')
\end{split}
\end{equation}
for any $ g \in Hom(\D, \B)$ such that $ g(x') = X_{Id}$, $g(z') = Z$, 
and $g(w') = W$; to do this we simply
compute the value of ${\sf c}_{(\Gamma)} g(s'\cdot t')$. This part of the proof is taken verbatim from Pigozzi \cite{P}. For the purpose of this
computation we make use of the following property of ordinals: if
$\Delta$ is any non-empty set of ordinals, then $ \bigcap \Delta$ is the
smallest ordinal in $ \Delta$, and if, in addition, $\Delta $ is
finite, then $\bigcup \Delta$ is the largest element ordinal in 
$\Delta$. Also, in this computation we shall assume that $\varphi$
always represents an arbitrary sequence in ${^{\alpha}{\alpha}}$.
Then, setting
$$ \Delta \varphi = \Gamma \sim \varphi[\Gamma \sim \{0, 1 \}]$$
for every $\varphi$, we successively compute:
$$ {\sf c}_1Z = \{ \varphi : |\Delta \varphi | = 2, \varphi_0 = \bigcap
\Delta \varphi \} \cap {\sf c}_{(\Gamma)} \{Id\},$$
\begin{equation*}
\begin{split}
& (X_{Id} \sim {\sf c}_1Z)
\cap {\sf c}_{(\Gamma)} \{ Id\} = \\
 & \{ \varphi : |\Delta \varphi | = 2, \varphi_0 = \bigcup \Delta
\varphi, \varphi_1 = \bigcap \Delta \varphi \}\cap {\sf c}_{(\Gamma)} \{Id\},
\end{split}
\end{equation*}
and, finally,
\begin{equation}\label{p11}
\begin{split}
{\sf c}_0(X_{Id}\sim {\sf c}_1Z)
\cap {\sf c}_{(\Gamma)} \{Id\} = \\
 & \{ \varphi : |\Delta \varphi | = 2, \varphi_1 =
 \bigcap \Delta \varphi \}\cap {\sf c}_{(\Gamma)} \{
Id\}.
\end{split}
\end{equation}
Similarly, we obtain
\begin{equation*}
\begin{split}
{\sf c}_0(X_{Id}\sim {\sf c}_1W)
\cap {\sf c}_{(\Gamma)} \{Id\} = \\
 & \{ \varphi : |\Delta \varphi | = 2, \varphi_1 =
 \bigcup \Delta \varphi \}\cap {\sf c}_{(\Gamma)} \{Id\}.
\end{split}
\end{equation*}
The last two formulas together give
\begin{equation}\label{p12}
\begin{split}
{\sf c}_0(X_{Id}\sim {\sf c}_1Z) \cap {\sf c}_0 (X_{Id} \sim {\sf c}_1W) \cap {\sf c}_{(\Gamma)} \{Id\} = 0.
\end{split}
\end{equation}
Continuing the computation we successively obtain:
$$ {\sf c}_1Z \cap {\sf d}_{01} = \{ \varphi : |\Delta \varphi | = 2,
\varphi_0 = \varphi_1 =
 \bigcap \Delta \varphi \}\cap {\sf c}_{(\Gamma)}\{Id\},$$
$$ {\sf s}^0_1{\sf c}_1Z = \{ \varphi : |\Delta \varphi | = 2,
\varphi_1 = \bigcap \Delta \varphi \}\cap {\sf c}_{(\Gamma)} \{ Id\},$$
$$ {\sf c}_1Z \cap {\sf s}^0_1{\sf c}_1Z = \{ \varphi : |\Delta \varphi | = 2,
\varphi_0 = \varphi_1 =
 \bigcap \Delta \varphi \}\cap {\sf c}_{(\Gamma)} \{Id\};$$
hence we finally get
\begin{equation}\label{p13}
\begin{split}
{\sf c}_0{\sf c}_1({\sf c}_1Z \cap {\sf s}^0_1{\sf c}_1Z \cap \sim {\sf d}_{01}) =  
{\sf c}_0{\sf c}_1 0 = 0,
\end{split}
\end{equation}
and similarly we get
\begin{equation}\label{p14}
\begin{split}
{\sf c}_0{\sf c}_1({\sf c}_1W \cap {\sf s}^0_1{\sf c}_1W \cap \sim {\sf d}_{01}) = 0.
\end{split}
\end{equation}
Now take $g$ to be any homomorphism from $\D$ into $\B$ such that
$g(x') = X_{Id}$, $g(z') = Z $ and $ g(w') = W$. Let $a=g(s'\cdot t')$. Then from the above 
$$ a\cap {\sf c}_{(\Gamma)} \{Id\} = \emptyset .$$
Then applying ${\sf c}_{(\Gamma)}$ to both sides of this equation we get
$${\sf c}_{(\Gamma)}a \cap {\sf c}_{(\Gamma)} \{Id\} = \emptyset . $$
Thus (\ref{p10}) holds.
Now there exists $ \Gamma \subseteq_{\kappa} \alpha$ 
and an interpolant $u' \in \D^{(x')},$  that is
$$r'\leq u'\leq {\sf c}_{(\Gamma)}(s'\cdot t').$$
There also exist $ Y, Z, W \subseteq
{^{\alpha}{\alpha}}$ such that (\ref{p6}) and (\ref{p7}) hold. Take
any $k \in Hom (\D, \B)$ such that $ k(x') = X_{Id}$, 
$ k(y') = Y$, $k(z') = Z$, and $k(w') = W.$
This is possible by the freeness of $\D.$ 
Then using the fact that $X_{Id} \cap k(r')$ is non-empty by (\ref{p6}) we get
$$ X_{Id} \cap k(u') = k(x'\cdot u') \supseteq k(x'\cdot r')
\neq 0.$$
And using the fact that $X_{Id} \sim
{\sf c}_{(\Gamma)}k(s'\cdot t')$ is non-empty by (\ref{p7}) we get
$$ X_{Id} \sim k(u') = k( x'\cdot -u') \supseteq
k(x'\cdot -{\sf c}_{(\Gamma)}(s'\cdot t')) \neq 0.$$
However, in view of (\ref{p4}), it is impossible for $X_{Id}$ 
 to intersect both $k(u')$ and its complement since
$ k(u') \in C$ and $X_{Id}$ is an atom; 
to see that $k(u')$ is indeed contained in $C$ recall that
$ u' \in \D^{(x')}$, and then observe that because of
(\ref{p5}) and the fact that $X_{Id} \in C$ we
must  have $k[\D^{(x')}] \subseteq C.$
This contradiction shows that $\K$ does not have the amalgamation property with respect to $PEA$. By this the proof is complete. 
\end{demo}
Other algebraic logics to which our proof  applies are Halmos' quasi-polyadic equality algebras  and Lucas' $\kappa$ extended cylindric algebras \cite{HMT2} p.267.
In particular, many varieties of those fail to have the amalgamation property.
We recall that Halmos quasi-polyadic algebras are of the form
$$\A=\langle A,+,\cdot ,-,0,1,{\sf c}_{(\Gamma)},{\sf s}_{\tau}, {\sf d}_{ij} \rangle_{i,j\in \alpha, \Gamma\subseteq_{\omega}\alpha ,\tau\in G_{\omega}}$$
while Lucas, $\kappa$ extended cylindric algebras are of the form
$$\A=\langle A,+,\cdot ,-,0,1,{\sf c}_{(\Gamma)}, {\sf d}_{\E} \rangle_{\Gamma\subseteq_{\kappa} \alpha, \E\in N}.$$
Both classes of (abstract) algebras are defined by a finite schema analogous to Halmos' schemas restricted to the appropriate similarity type, cf. Def 1.1.
The representable algebras are defined as subdirect product of set algebras. In those two cases the class of representable algebras, 
as opposed to the class of abstract algebras, is not finite 
schema axiomatizable. The methods of Andreka  in \cite{A} can be used to prove this (the proof though is not trivial). But in those two cases the class of representable algebras forms a variety
and using our proof it can be easily shown that any variety containing the representable algebras such that its cylindric reduct 
satisfies the cylindric axioms 
fails to have the amalgamation property.  In particular, in both of these cases, both the variety of abstract algebras as well as that of the 
representable algebras fail 
to have the amalgamation property.

\end{document}